\documentclass[12pt, a4paper]{article}

\usepackage{amsfonts, amsmath, amssymb, amsthm}
\usepackage{changepage}
\usepackage{enumitem}
\usepackage{fullpage}
\usepackage{hyperref}
\usepackage{mathdots}
\usepackage{parskip}
\usepackage{tikz}
\usepackage{xurl} 
\usepackage{lipsum}
\usepackage{comment}

\usetikzlibrary{arrows}

\DeclareMathOperator{\sgn}{sgn}

\theoremstyle{plain}
\newtheorem{theorem}{Theorem}

\newtheorem{corollary}[theorem]{Corollary}
\newtheorem{proposition}[theorem]{Proposition}
\newtheorem{lemma}[theorem]{Lemma}

\begin{document}

\title{On other two representations of the C-recursive integer sequences by terms in modular arithmetic}

\author{Mihai Prunescu \footnote{Research Center for Logic, Optimization and Security (LOS), Faculty of Mathematics and Computer Science, University of Bucharest, Academiei 14, Bucharest (RO-010014), Romania; and Simion Stoilow Institute of Mathematics of the Romanian Academy, Research unit 5, P. O. Box 1-764, Bucharest (RO-014700), Romania. E-mail: {\tt mihai.prunescu@imar.ro}, {\tt mihai.prunescu@gmail.com}.  Partially supported by Bitdefender (Research in Pairs in Bucharest Program).}}
\date{}
\maketitle

\begin{abstract} An integer sequence that is defined  by initial values and a linear recurrence with constant integer coefficients, can be represented by the difference of two arithmetic terms containing exponentiation. All constants occuring in the term are integers. While in the paper {\it On the representation of C-recursive integer
sequences by arithmetic terms} by Prunescu and Sauras-Altuzarra, the terms consist of the remainder operation, applied on a division; the representations shown here are a division applied to a remainder operation, respectively the composition of two remainder operations.   \end{abstract}

AMS Subject Classification: 11B37 (primary), 11Y55, 39A06.

Keywords: C-recursive sequence, modular arithmetic, term representation, Pell's equation, Lucas sequence, Narayana's Cows.

\section{Introduction}\label{Introduction}

Given a ring $ R $, a sequence of terms of $ R $ that satisfies a \textit{homogeneous linear recurrence} of constant \textit{coefficients} in $ R $ (i.e.\ a sequence $ s \in R^\mathbb{N} $ such that, for some positive integer $ d $ and $ d $ elements $ \alpha_1 $, $ \dots $, $ \alpha_d $ of $ R $, we have that $ \alpha_d $ is non-zero and $$ s ( n + d ) + \alpha_1 s ( n + d - 1 ) + \dots + \alpha_d s ( n ) = 0 $$ for every integer $ n \geq 0 $) is said to be \textbf{C-recursive} of \textbf{degree} $ d $ (cf.\ Petkovšek \& Zakrajšek \cite[Definition 1]{PetkovsekZakrajsek}).

In a recent paper, Prunescu and Sauras-Altuzarra \cite{PrunescuSaurasCRec1} shown that a C-recursive sequences of natural numbers $s(n)$ always satisfies a relation of the form:
$$\forall\, n \geq n_0\,\,\,\,s(n) = \left \lfloor \frac{A(b,n)}{B(b,n)} \right \rfloor \bmod b^n.$$
Here $A(b,n)$ and $B(b,n)$ are exponential polynomials such that $n$ occurs only in some exponents of the basis $b \in \mathbb N$, those exponents being polynomials of $n$. The representations work for sufficiently large $b$, which depends only on the sequence. We will call the above representation, the {\bf quotient - remainder} representation.

In this paper I will show that other two representations work:
$$ \forall\,n \geq n_0\,\,\,\, s(n) = \left \lfloor \frac{A(b,n) \bmod B(b,n)}{b^{(d-1)n}}\right \rfloor,$$
which will be called the {\bf remainder - quotient} representation, respectively:
$$ \forall\, n \geq n_0\,\,\,\, s(n) = A(b,n) \bmod B(b,n) \bmod b^n,$$
which will be called the {\bf remainder - remainder} representation. It must be specified that the expression:
$$A \bmod B \bmod C$$
will be always understood as:
$$(A \bmod B) \bmod C.$$

It is worth to mention here that the representation presented in \cite{PrunescuSaurasCRec1} has been seen by the authors as a step forward in the framework of Mazzanti and Marchenkov. Our research program deals with representing concrete functions which are Kalmar elementary (functions which can be computed in deterministic time bounded by an exponential tower of fixed height, in the length of the input) as arithmetic terms. So the operations used to compose the arithmetic term were subject of the restrictions of this framework. In particular, one may use only the {\bf modified} subtraction $x \dot - y$, which equals $x-y$ only if $x \geq y$ and is $0$ otherwise.

In the present paper no modified subtraction is used. The operations used are the addition $x + y$, the (full) subtraction $x - y$, the product $xy$, the quotient and the remainder $x = qy + m$ denoted by $q = \lfloor x/y \rfloor$ and by $m = x \bmod y$ and the exponential function $b^x$, all defined over the set of integer numbers $\mathbb Z$. It is of crucial importance that the remainder $m = x \bmod y$ satisfies $0 \leq m < |y|$, as it is also of crucial importance that in some expressions $A(b,n) \bmod B(b,n)$ the value of $A(b,n)$ will always be a negative integer, while the value of this expression will be always positive. The only one restriction concerns expressions of the form:
$$b^{f(n),}$$
where $b$ and $f(n)$ should always be natural numbers. 

It might be of interest how this paper arose. Well, it begun with a Romanian Easter Story. On May 5-th 2024, the day of the Romanian Easter this year, shortly before the paper \cite{PrunescuSaurasCRec1} was submitted to arXiv, the I wrote the following post in a Facebook group called {\it Arxiv Math} (whose members are some Romanian mathematicians) the following question: {\it On the Easter Day, problems with Fibonacci's sequence are very recommended, as this sequence works with bunnies. Please compute the integer part of the following rational number:}
$$\frac{10^{110}}{10^{20} - 10^{10} - 1}. $$
The result of this computation is the number:
$$\textbf{1}000000000\textbf{1}000000000\textbf{2}000000000\textbf{3}000000000\textbf{5}000000000\textbf{8} \dots $$ $$ \dots 00000000\textbf{13}00000000\textbf{21}00000000\textbf{34}00000000\textbf{55}  $$
which encodes a starting segment of the Fibonacci sequence, see Prunescu and Sauras-Altuzarra \cite{PrunescuSaurasCRec1}. Professor Adrian Atanasiu commented to my post the same day. He wrote: {\it I denoted the number $10^{10}$ with $a$ and I computed the polynomial division $a^{11}/ (a^2 - a - 1)$. I don't want to write down the quotient, but the remainder is} $$\textbf{89}00000000\textbf{55}$$ {\it and contains Fibonacci numbers as well.} The second answer to my post belongs to Professor Marcel \c{T}ena, who observed that the following polynomial identity in Fibonacci numbers holds for all $n \geq 2$: 
$$ a^n=(a^2-a-1)(F(0)a^{n-2}+F(1)a^{n-3}+\dots +F(n-3)a+F(n-2))+$$ $$+ F(n-1)a+F(n-2).$$
This is a particular case of the general identity that you will find in Section \ref{section2}. It was clear that not only the quotient given by the generating function of a C-recursive sequence encodes a segment of the sequence, but that the remainder encodes a segment of the sequence as well. In this paper I show how to extract the information contained in the remainder, for any C-recursive sequence of natural numbers. I also show how to apply this result to represent arbitrary C-recurrent sequences of integers, including the mysterious oscillating sequences, which contains infinitely many alternations of sign. 

While in the paper Prunescu and Sauras-Altuzarra \cite{PrunescuSaurasCRec1} some Analysis was used, and arguments went frequently about the generating function, its radius of convergence in $0$, the value of a truncated series, and so on; the tools used in the present paper are only Algebra of polynomials and some elementary Number Theory. 

As always, the story does not end here. One week before the present paper was submitted to arXiv, Professor Michael Wallner informed my co-author from \cite{PrunescuSaurasCRec1}, Professor Sauras-Altuzarra, about a very interesting paper by Joseph M. Shunia, called {\it Polynomial Quotient Rings and Kronecker Substitution for Deriving Combinatorial Identities}, see \cite{Shunia}. Among a plethora of interesting ideas and results, Joseph Shunia also proves  a remainder - remainder representation for C-recursive sequences. Surprisingly, both the construction and the result are different from the remainder - remainder representation done here.   From Shunia's paper, I learnt that a quotient - remainder representation of the Fibonacci sequence was already presented 2016 in a blog post by Paul Hankin, see \cite{Hankin}, and a remainder - remainder representation for Fibonacci was intuitively found (see the same blog post) by the contributor Far{\'e} Rideau. In his post, Hankin asked if such representations are possible for other sequences. The paper \cite{PrunescuSaurasCRec1}, Joseph Shunia's paper \cite{Shunia} and the present paper contain different positive answers to this question.

A final comment on how to possibly apply these results. I believe that the most practical representation is the remainder - remainder representation. Not only that one can use routines of modular arithmetic to fast compute $A(b,n) \bmod B(b,n)$, but because the division with $b^{n}$, one has to compute only the last $n$ base $b$-digits of the final computation.

\section{Polynomial identities}\label{section2}

Let $R$ be some ring and $s \in R^{\mathbb N}$ be a C-recursive sequence of degree $d$. We associate with the sequence $s$ and with its recurrence rule 
$$ s ( n + d ) + \alpha_1 s ( n + d - 1 ) + \dots + \alpha_d s ( n ) = 0 $$
with $\alpha_1, \dots, \alpha_d \in R$, $\alpha_d \neq 0$, the following polynomials in $R[X]$. Let $B_s(X) \in R[X]$ be the polynomial defined as:
$$B_s(X) = \alpha_d + \alpha_{d-1} X + \dots + \alpha_1 X^{d-1} + X^d.$$
For every $k \in \mathbb N$ let $S_k(X)$ be the polynomial defined as:
$$S_k(X) = s(k) + s(k-1) X + \dots + s(0) X^k. $$

\begin{theorem}\label{decomposition}
    Let $s \in R^{\mathbb N}$ be a C-recursive function of degree $d \geq 2$. Then there is a unique polynomial $A(X) \in R[X]$ such that for every natural number $k \geq d-1$, there is a unique polynomial $C_k(X) \in R[X]$ such that:
    \begin{enumerate}
        \item The following polynomial identity is true:
        $$ B(X) S_k(X) = C_k(X) + X^{k + 1} A(X).$$ 
        \item $\deg(A(X)) < d$. Moreover, if $A(X) = a_{d-1} + a_{d-2}X + \dots + a_0X^{d-1}$, then $a_0 = s(0)$, $a_1$ is a fixed linear combination of $s(0)$ and $s(1)$, and in general $a_i$ is a fixed linear combination of  $s(0)$, $s(1)$, $\dots$, $s(i)$ for all $i = 0, \dots, d-1$. 
        \item $\deg(C_k(X)) < d$. If $C_k(X) = c_0 + c_1X + \dots + c_{d-1} X^{d-1}$ then $c_0 = \alpha_ds(k)$, $c_1$ is a fixed linear combination of $s(k)$ and $s(k-1)$, and in general $c_i$ is a fixed linear combination of  $s(k)$, $s(k-1)$, $\dots$, $s(k-i)$ for all $i = 0, \dots, d-2$. In these linear combinations, the coefficients  do not depend on $k$. 
        \item The coefficient $c_{d-1} = - s(k+1)$. 
    \end{enumerate}
\end{theorem} 

{\bf Proof} Most conclusions are immediate. The coefficient of $X^{d-1}$ in $B(X)S_k(X)$ is 
$$\alpha_1 s(k) + \alpha_2 s(k-1) + \dots + \alpha_d s(k - d + 1) = - s(k+1)$$
according to the recurrence relation. Also, if $d \leq u \leq k$ then the coefficient of $X^u$ in the product $B(X)S_k(X)$ will be: 
$$ s(k - u + d) + \alpha_1 s(k - u + d - 1) + \dots + \alpha_d s(k - u) = 0$$
by the same relation of recurrence. \qed 

We observe that the polynomials given in Theorem \ref{decomposition} may be defined as: 
\begin{eqnarray*}
    C_k(X) &=& B_s(X) S_k(X) \bmod X^d, \\
    A_s(X) &=& B_s(X) S_{d-1}(X) \,\, / \,\, X^d.
\end{eqnarray*}

\section{Inequalities} 

The following Lemma was proved in Prunescu \& Sauras-Altuzarra, \cite{PrunescuSaurasCRec1}. 

\begin{lemma}\label{LemmaGeneralInequality} If $ s \in \mathbb{C}^{\mathbb{N}} $ is C-recursive, then there is an integer $ c \geq 1 $ such that $ | s(n) | < c^{n+1} $ for every integer $ n \geq 0 $. \end{lemma} 

The following Lemma is similar with one proved in the same paper, but is given here slightly more generally: 

\begin{lemma}\label{Lemman/s} For  $ b , c, s \in \mathbb R$ such that $s \geq 1$ and  $ b > c^{2s} > 1 $, one has that $ c^{n + 1} < b^{n / s} $ for every integer $ n \geq 1 $. \end{lemma} 

{\bf Proof} As $\frac{n}{2s} + \frac{1}{2s} \leq \frac{n}{s}$ is equivalent with $n \geq 1$, one has:
$$c^{n+1} < {\left (\sqrt[2s]{b} \right )}^{n+1} = b^{\frac{n}{2s} + \frac{1}{2s} } \leq b^{n/s}.$$
\qed 

\begin{lemma}\label{Lemman/3}
    If $ s \in \mathbb{C}^{\mathbb{N}} $ is C-recursive, then there is an integer $ b \geq 1 $ such that $ | s(n) | < b^{n/3} $ for every integer $ n \geq 0 $.
\end{lemma} 

{\bf Proof} From Lemma \ref{LemmaGeneralInequality}, there is a $ c \in \mathbb N $ such that for all $n \geq 1$, $ | s(n) | < c^{n+1} $. From Lemma \ref{Lemman/s} applied for $s = 3$ there is a $b \in \mathbb N$ such that for all $n \geq 1$, $c^{n+1} < b^{s/3}$. The conclusion follows.

Ne we apply them in the context of Theorem \ref{decomposition}. 

\begin{lemma}\label{Lemmacoeffinequalities}
Let $s \in \mathbb N ^ {\mathbb N}$ be a recursive sequence of degree $d \geq 2$, the recurrence relation having coefficients in $\mathbb Z$. For every natural number $k \geq d-1$, let
$$C_k(X) = c_0(k) + c_1(k) X + \dots + c_{d-1}(k) X^{d-1}$$
be the polynomial such that
$$ C_k(X) = B_s(X)S_k(X)  \bmod X^{d}. $$
Then there exists $b \in \mathbb N$ such that
$$|c_i(k)| < b^{k/3}$$
for all $i \in \{0, 1, \dots, d-1\}$ and for all $k \geq d-1$.
\end{lemma}

{\bf Proof} $C_k(X)$ is the polynomial whose existence is stated in Theorem \ref{decomposition}. By the same theorem, for all $i \in \{0, 1, \dots, d-1\}$, $c_i$ is a fixed linear combination of the numbers $s(k), s(k-1), \dots s(k-i)$ with coefficients in $\mathbb Z$, say:
$$c_i(k) = A_{i0} s(k) + A_{i1} s(k-1) + \dots + A_{i i} s(k-i).$$

According to the Lemma \ref{Lemman/3}, as every sequence $(s(n-i)\,|\, n - i \geq 0)$ is C-recursive,  there is $b_1 \in \mathbb N$ sufficiently large such that for all $k \geq d-1$ and for all $i \in \{0, 1, \dots, d-1\}$,
$$|s(k-i)| < b_1^{k/3}.$$
Replacing this value of $b_1$ with a bigger one, if necessary, we can assure that for all $i \in \{0, 1, \dots, d-1\}$,
$$|A_{i0}| + |A_{i1}| + \dots + |A_{i i}| < b_1^{1/3}.$$
It follows that for all $i \in \{0, 1, \dots, d-1\}$,
$$|c_i(k)| < b_1^{k/3} b_1^{1/3} = b_1^{(k+1)/3}.$$
According to Lemma \ref{Lemman/s} for $s = 1$, $c = b_1^{1/3}$ and $b > b_1^{2/3}$, for all $k \geq 1$
$$b_1^{(k+1)/3} < b^{k/3}.$$
So for all  $k \geq d-1$ and for all $i \in \{0, 1, \dots, d-1\}$,
$$|c_i(k)| < b^{k/3}.$$
\qed 

Consider a finite sequence of positive real numbers
$$r_s > r_{s-1} > \dots > r_1 > 0. $$
For arbitrary coefficients $c_0, c_1, \dots c_{s-1} \in \mathbb R$, the {\bf generalized monic polynomial} with exponents $r_i$ and coefficients $c_i$ is the function:
$$f(Y) = Y^{r_s} + c_{s-1}Y^{r_{s-1}} + \dots + c_1 Y^{r_1} + c_0.$$
We say that $\deg f(Y) = r_s$.

\begin{lemma}\label{LemmaGeneralizedPolynomials}
    Let $f(Y)$ and $g(Y)$ two generalized monic polynomials with $$\deg f(Y) > \deg g(Y),$$ and with not necessarily the same number of summands. Then there is some $y_0 \in \mathbb R$ such that for all $y > y_0$,
    $$ f(y) > g(y) > 0.$$
\end{lemma} 

{\bf Proof} We observe that $\lim_{y \rightarrow \infty} g(y) = + \infty$ and that $\lim_{y \rightarrow \infty} (f(y) - g(y))  = + \infty$. The conclusion follows. \qed 

\section{The remainder -- quotient representation}

\begin{theorem}\label{Theoremmoddivrep}
    Let $s \in \mathbb N^{\mathbb N}$ be a C-recursive sequence of degree $d \geq 2$. If $B(X) = B_s(X)$ is the polynomial associated to its recurrence relation and $A(X) = A_s(X)= B_s(X)S_{d-1}(X) / X^d$, then there is $n_0$ such that for a sufficiently large $b \in \mathbb N$:
    $$\forall \, n\geq n_0\,\,\,\, s(n) = \left \lfloor 
    \frac{\left ( b^{n(d-2) + \left \lceil n/2 \right \rceil} + b^{n^2}A(b^n ) \right ) \bmod B(b^n)}{b^{n(d-1)}} \right \rfloor.$$
\end{theorem}

{\bf Proof} Recall the relation from Theorem \ref{decomposition} for $k = n-1$:
$$ B(X) S_{n-1}(X) = C_{n-1}(X) + X^{n} A(X).$$ 
It follows that:
$$-C_{n-1}(X) \equiv X^n A(X) \mod B(X).$$
If we develop $C_{n-1}(X)$, we get:
$$s(n) X^{d-1} - c_{d-2}(n-1) X^{d-2} - \dots - c_0(n-1) \equiv X^n A(X) \mod B(X).$$
We will make the substitution $X = b^n$ for a sufficiently large $b \in \mathbb N$. We consider the following generalized monic polynomials:
\begin{eqnarray*}
f(Y) &=& Y^{d-1 + \frac{1}{3}} + Y^{d - 2 + \frac{2}{3}} + Y^{{d-2} + \frac{1}{3}} + \dots + Y^{\frac{1}{3}} \\
g(Y) &=& Y^{d-1} - Y^{d - 2 + \frac{2}{3}} - Y^{{d-2} + \frac{1}{3}} - \dots - Y^{\frac{1}{3}} \\
h(Y) &= & Y^{d - 2 + \frac{1}{2}} - Y^{{d-2} + \frac{1}{3}} - \dots - Y^{\frac{1}{3}} 
\end{eqnarray*}
As $B(Y)$ is a monic polynomial of degree $d$, according to Lemma \ref{LemmaGeneralizedPolynomials} there is some $y_0 \in \mathbb R$ such that for all $y > y_0$, $B(y) > f(y)$. 
Also, there are some $y_1, y_2 \in \mathbb R$ such that for all $y > y_1$, $g(y) > 0$ and for all $y > y_2$, $h(Y) > 0$.

We choose $b \in \mathbb N$ such that according to Lemma \ref{Lemmacoeffinequalities}, $|c_i(n)| < b^{n/3}$ for all $n \geq 1$ and for all $i \in \{0, 1, \dots, d-1\}$. Also, by eventually replacing $b$ with a bigger one, $b > y_i$ for $i = 0, 1, 2$. We see that: 
$$s(n) b^{n(d-1)} - c_{d-2}(n-1) b^{n(d-2)} - \dots - c_0(n-1) \equiv b^{n^2} A(b^n) \mod B(b^n).$$
We add $b^{n(d-2) + \left \lceil n/2 \right \rceil}$ to this relation. It follows: 
$$s(n) b^{n(d-1)} + b^{n(d-2) + \left \lceil n/2 \right \rceil} - c_{d-2}(n-1) X^{n(d-2)} - \dots - c_0(n-1) \equiv $$ $$ \equiv (b^{n(d-2) + n/2} + b^{n^2} A(b^n)) \mod B(b^n).$$
We observe that: 
$$s(n) b^{n(d-1)} + b^{n(d-2) + \left \lceil n/2 \right \rceil} - c_{d-2}(n-1) X^{n(d-2)} - \dots - c_0(n-1) < $$
$$< b^{n(d-1) + n / 3} + b^{n(d-2) + 2n/3}  + b^{n(d-2) + n/3} + \dots + b^{n/3} = f(b^n) < B(b^n). $$
For the other inequality, we must remember that some terms $s(n)$ with $n \geq d-1$ might be equal $0$. So: 
$$s(n) b^{n(d-1)} + b^{n(d-2) + \left \lceil n/2 \right \rceil} - c_{d-2}(n-1) X^{n(d-2)} - \dots - c_0(n-1) > $$
$$ > b^{n(d-2) +  n/2 } - b^{n(d-2) + n/3} - \dots - b^{n/3} = h(b^n) > 0. $$
We conclude that
$$\left (b^{n(d-2) + \left \lceil n/2 \right \rceil} + b^{n^2} A(b^n) \right ) \bmod B(b^n) =  $$
$$= s(n) b^{n(d-1)} + b^{n(d-2) + \left \lceil n/2 \right \rceil} - c_{d-2}(n-1) X^{n(d-2)} - \dots - c_0(n-1).$$
We observe that: 
$$ b^{n(d-2) + \left \lceil n/2 \right \rceil} - c_{d-2}(n-1) X^{n(d-2)} - \dots - c_0(n-1) > $$
$$b^{n(d-2) + n/2 } - b^{n(d-2) + n/3} - \dots - b^{n/3} = h(b^n) > 0,$$
and also that:
$$ b^{n(d-1)} - \left ( b^{n(d-2) + \left \lceil n/2 \right \rceil} - c_{d-2}(n-1) X^{n(d-2)} - \dots - c_0(n-1)  \right ) > $$
$$ > b^{n(d-1)} - b^{n(d-2) +  2n/3}  - b^{n(d-2) + n/3} - \dots - b^{n/3} = g(b^n) > 0.$$
We conclude that:
$$0 < b^{n(d-2) + \left \lceil n/2 \right \rceil} - c_{d-2}(n-1) X^{n(d-2)} - \dots - c_0(n-1) < b^{(d-1)n},$$
so that:
$$\left \lfloor \frac{(b^{n(d-2) + \left \lceil n/2 \right \rceil} + b^{n^2} A(b^n)) \bmod B(b^n)}{b^{(d-1)n}} \right \rfloor = s(n).$$
\qed 

\section{The remainder -- remainder representation}

\begin{theorem}\label{Theoremmodmodrep}
    Let $s \in \mathbb N^{\mathbb N}$ be a C-recursive sequence of degree $d \geq 2$. If $B(X) = B_s(X)$ is the polynomial associated to its recurrence relation and $A(X) = A_s(X)= B_s(X)S_{d-1}(X) / X^d$, then there is $n_0$ such that for a sufficiently large $b \in \mathbb N$:
    $$\forall \, n\geq n_0\,\,\,\, s(n) = \frac{1}{|\alpha_{d} |}  \left (b^{n(d-1) + \left \lceil n/2 \right \rceil} -  b^{n(n+1)} \sgn(\alpha_d) A(b^n) \right )\bmod B(b^n)       \bmod b^n.$$
\end{theorem}

{\bf Proof} Recall the relation from Theorem \ref{decomposition} for $k = n$:
$$ B(X) S_{n}(X) = C_{n}(X) + X^{n+1} A(X).$$ 
It follows that:
$$C_{n}(X) \equiv - X^{n+1} A(X) \mod B(X).$$
If we develop $C_{n}(X)$, we get:
$$\alpha_d s(n) + c_1(n) X + \dots + c_{d-1}(n) X^{d-1} \equiv - X^{n+1} A(X) \mod B(X).$$ 
Let $s = \alpha_d / |\alpha_d| \in \{-1, 1\}$. The whole relation is multiplied by $s$. 
$$|\alpha_d | s(n)  + sc_1(n) X + \dots + sc_{d-1}(n) X^{d-1} \equiv - s X^{n+1} A(X) \mod B(X).$$ 
We will make the substitution $X = b^n$ for a sufficiently large $b \in \mathbb N$. We consider the following generalized monic polynomials:
\begin{eqnarray*}
f(Y) &=& Y^{d-1 + \frac{2}{3}} + Y^{d - 1 + \frac{1}{3}} + Y^{{d-2} + \frac{1}{3}} + \dots + Y^{\frac{1}{3}} \\
g(Y) &=& Y^{d-1 + \frac{1}{2}} - Y^{d - 1 + \frac{1}{3}} - Y^{{d-2} + \frac{1}{3}} - \dots - Y^{\frac{1}{3}} \\
h(Y) &= & Y^{d - 1 + \frac{1}{2}} - Y^{{d-1} + \frac{1}{3}} - \dots - Y^{\frac{4}{3}} 
\end{eqnarray*} 
As $B(Y)$ is a monic polynomial of degree $d$, according to Lemma \ref{LemmaGeneralizedPolynomials} there is some $y_0 \in \mathbb R$ such that for all $y > y_0$, $B(y) > f(y)$. 
Also, there are some $y_1, y_2 \in \mathbb R$ such that for all $y > y_1$, $g(y) > 0$ and for all $y > y_2$, $h(Y) > 0$.

We choose $b \in \mathbb N$ such that according to Lemma \ref{Lemmacoeffinequalities}, $|c_i(n)| < b^{n/3}$ for all $n \geq 1$ and for all $i \in \{0, 1, \dots, d-1\}$. Also, by eventually replacing $b$ with a bigger one, $b > y_i$ for $i = 0, 1, 2$. We see that: 
$$|\alpha_d | s(n)  + sc_1(n) b^n + \dots + sc_{d-1}(n) b^{n(d-1)}  \equiv - s b^{n(n+1)} A(b^n) \mod B(b^n).$$
We add $b^{n(d-1) + \left \lceil n/2 \right \rceil}$ to this relation. It follows: 
$$|\alpha_d | s(n)  + sc_1(n) b^n + \dots + sc_{d-1}(n) b^{n(d-1)} 
+ b^{n(d-1) + \left \lceil n/2 \right \rceil} \equiv $$$$ \equiv
b^{n(d-1) + \left \lceil n/2 \right \rceil} - s b^{n(n+1)} A(b^n) \mod B(b^n).$$
We observe that
$$|\alpha_d | s(n)  + sc_1(n) b^n + \dots + sc_{d-1}(n) b^{n(d-1)} 
+ b^{n(d-1) + \left \lceil n/2 \right \rceil} < $$
$$ < b^{n/3} + b^{n(1 + 1/3)} + \dots + b^{n(d - 1 + 1/3)} + b^{n(d-1 + 2/3)} < f(b^n) < B(b^n),$$
$$|\alpha_d | s(n)  + sc_1(n) b^n + \dots + sc_{d-1}(n) b^{n(d-1)} 
+ b^{n(d-1) + \left \lceil n/2 \right \rceil} > $$
$$> - b^{n/3} - b^{n(1 + 1/3)} - \dots - b^{n(d - 1 + 1/3)} + b^{n(d-1 + 1/2)} > g(b^n) > 0. $$
We conclude that
$$\left (b^{n(d-1) + \left \lceil n/2 \right \rceil} - s b^{n(n+1)} A(b^n) \right )\bmod B(b^n) = $$
$$ = |\alpha_d | s(n)  + sc_1(n) b^n + \dots + sc_{d-1}(n) b^{n(d-1)} + b^{n(d-1) + \left \lceil n/2 \right \rceil}.$$
Also,
$$  sc_1(n) b^n + \dots + sc_{d-1}(n) b^{n(d-1)} + b^{n(d-1) + \left \lceil n/2 \right \rceil} > $$
$$ > - b^{n(1 + 1/3)} - \dots - b^{n(d-1+1/3)} + b^{n(d-1+1/2)} > 0.$$
As this number is also a multiple of $b^n$, we conclude that:
$$|\alpha_d | s(n) = \left ( \left (b^{n(d-1) + \left \lceil n/2 \right \rceil} - s b^{n(n+1)} A(b^n) \right )\bmod B(b^n)      \right ) \bmod b^n.$$
\qed

\section{Sequences of integers}

\begin{lemma}\label{Lemmarecurrencesumisprod}
    Let $R$ be a ring. Let $s, t \in R^{\mathbb N}$ be two sequences, that satisfy the following relations of recurrence for all $n \geq 0$:
    \begin{eqnarray*}
       \sum _{a=0}^d \alpha_a s(n+a) &=& 0, \\
        \sum _{b=0}^g \beta_b t(n+b) &=& 0, 
    \end{eqnarray*}
    where $\alpha_i, \beta_j \in R$ and $d, g \geq 1$. Consider the polynomials $B_s(X) = \alpha_d X^d + \alpha_{d-1} X^{d+n-1} + \dots + \alpha_0$ and $B_t(X) = \beta_gX^g + \beta_{g-1}X^{g-1} + \dots + \beta_0 $. Let
    $$B_{s+t}(X) := B_s(X) B_t(X) = \sum _{i=0}^{d+g} \gamma_i X^i.$$
    For all $n \in \mathbb N$ we define $u(n) = s(n) + t(n)$. Then the sequence $u$ satisfies the recurrence rule:
    $$  \sum _{i=0}^{d+g} \gamma_i u(n+i) = 0. $$
\end{lemma}

{\bf Proof} We denote the coefficients of $R_{s+t}$ with $\gamma_i$ as follows:
$$B_{s+t}(X) = \sum _{i=0}^{d+g} \left ( \sum _{a + b = i} \alpha_a \beta_b \right ) X^i = \sum _{i=0}^{d+g} \gamma_i X^i. $$
For some $n \in \mathbb N$ we compute:
$$  \sum _{i=0}^{d+g} \gamma_i u(n+i) =  \sum _{i=0}^{d+g} \gamma_i s(n+i) +  \sum _{i=0}^{d+g} \gamma_i t(n+i) = $$
$$ \sum _{b=0}^g \beta_b \sum _{a = 0}^d \alpha_a s(n+b+a) + 
\sum _{a=0}^d \alpha_a \sum _{b = 0}^g \beta_b t(n+a+b) = $$
$$ \sum _{b=0}^g \beta_b \cdot 0 + \sum _{a=0}^d \alpha_a \cdot 0 = 0. $$
\qed

\begin{theorem}\label{TheoremIntegers}
    For every C-recursive sequence $s \in \mathbb Z^{\mathbb N}$ there exists $c \in \mathbb N$ such that the C-recursive sequence defined as $u(n) = s(n) + c^{n+1}$ consists of natural numbers. 
\end{theorem} 

{\bf Proof} According to Lemma \ref{LemmaGeneralInequality}, there is an integer $c \geq 1$ such that for all $n \in \mathbb N$, $|s(n)| < c^{n+1}$. In concrete situations it might be of interest how one finds such numbers. From the proof of the Lemma (see \cite{PrunescuSaurasCRec1}) it is necessary and sufficient that:  $ d \left (|\alpha_1| + \dots + |\alpha_d| \right ) < c $, $|s(0)| < c$, $|s(1)| < c^2$, $\dots$, $|s(d-1)| < c^d$, where $\alpha_1, \dots, \alpha_d$ are such that $ s(n+d) =  - \alpha_1 s ( n + d - 1 ) - \dots - \alpha_d s ( n ) $ for every integer $ n \geq 0 $.

Further we observe that the sequence $u(n) := s(n) + c^{n+1}$ is indeed C-recursive according to Lemma \ref{Lemmarecurrencesumisprod}. Concretely, we take $B_s(X) = X^d + \alpha_{d-1} X^{d-1} + \dots + \alpha_0$ and $B_t(X) = X - c$. The recurrence relation of $u$ is given by the product $(X-c)B_s(X)$. The first $d+1$ terms $u(0) = s(0) + c$, $\dots$, $u(d) = s(d) + c^{d+1}$ can be computed explicitly. We also see that $u(n) = s(n) + c^{n+1} \geq 0$ for all $n \geq 0$. \qed 

In the following statements let $$U_n(X) := u(n) + u(n-1) X + \dots + u(0) X^n.$$ 
Also, $\gamma_{d+1}$ is the free term in the recurrence polynomial of the sequence $u(n)$. 

\begin{corollary}\label{Cormoddivforint}
     Let $s \in \mathbb Z^{\mathbb N}$ be a C-recursive sequence of degree $d \geq 2$. If $c > 0$ is sufficiently large, such that for all $n \in \mathbb N$, $|s(n)| < c^{n+1}$, $u(n) = s(n) + c^{n+1}$,  $B(X) = B_u(X)$ is the polynomial associated to the recurrence relation of the sequence $u$ and $A(X) = A_u(X)= B_u(X)U_{d}(X) / X^{d+1}$, then there is $n_0$ such that for sufficiently large $b \in \mathbb N$:
    $$\forall \, n\geq n_0\,\,\,\, s(n) = \left \lfloor 
    \frac{\left ( b^{n(d-1) + \left \lceil n/2 \right \rceil} + b^{n^2}A(b^n ) \right ) \bmod B(b^n)}{b^{nd}} \right \rfloor - c^{n+1}.$$
\end{corollary} 

\begin{corollary}\label{Cormodmodforint}
     Let $s \in \mathbb Z^{\mathbb N}$ be a C-recursive sequence of degree $d \geq 2$. If $c > 0$ is sufficiently large, such that for all $n \in \mathbb N$, $|s(n)| < c^{n+1}$, $u(n) = s(n) + c^{n+1}$,  $B(X) = B_u(X)$ is the polynomial associated to the recurrence relation of the sequence $u$ and $A(X) = A_u(X)= B_u(X)U_{d}(X) / X^{d+1}$, then there is $n_0$ such that for sufficiently large $b \in \mathbb N$:
      $$\forall \, n\geq n_0$$ $$ s(n) = \frac{1}{|\gamma_{d +1} |} \left (b^{nd + \left \lceil n/2 \right \rceil} -  b^{n(n+1)} \sgn(\gamma_{d+1}) A(b^n) \right )\bmod B(b^n)    \bmod b^n - c^{n+1}$$
\end{corollary}

\section{Sequences of low degree}

The goal of this section is two-fold. First we give explicit representation formulas for recurrent sequences of natural numbers of degree $2$ and $3$, respectively for recurrent sequences of integers of degree $2$. Second, we point out that in some particular cases, alternative representations work. 

\subsection{Degree 2, positive terms}

Let $s \in \mathbb N ^ {\mathbb N}$ be a sequence with initial terms $s(0)$ and $s(1)$ and with recurrence polynomial $B(X) = \alpha_2 + \alpha_1X + X^2 \in \mathbb Z[X]$. As the polynomial $S_n(X)$ ends in
$$ \dots + s(1) X^{n-1} + s(0) X^n,$$
we conclude that $$A(X) = s(0) X + (\alpha_1 s(0) + s(1)).  $$ According to Theorems \ref{Theoremmoddivrep} and \ref{Theoremmodmodrep}, the sequence $s$ has the following representations: 

$$\forall \, n\geq n_0\,\,\,\, s(n) = \left \lfloor 
    \frac{\left ( b^{\left \lceil n/2 \right \rceil} +b^{n^2}A(b^n) \right ) \bmod B(b^n)}{b^{n}} \right \rfloor,$$

$$\forall \, n\geq n_0\,\,\,\, s(n) = \frac{1}{|\alpha_{2} |}  \left (b^{n + \left \lceil n/2 \right \rceil} -  b^{n(n+1)} \sgn(\alpha_2) A(b^n) \right )\bmod B(b^n)       \bmod b^n.$$

 We call these representations {\it general}. Now we examine  the proofs of the Theorems \ref{Theoremmoddivrep} and \ref{Theoremmodmodrep} in order to find possible short-cuts for  this particular case. 

 For the remainder -- quotient representation, consider the identity
 $$ B(X) S_{n-1}(X) = C_{n-1}(X) + X^{n} A(X)$$ 
and observe that 
$$C_{n-1}(X) = \alpha_2 s(n-1) - X s(n).$$
There are no other terms between these. It follows that if $\alpha_2 < 0$ then 
$$(-\alpha_2) s(n-1) + X s(n) \equiv X^n A(X) \mod B(X),$$
$$(-\alpha_2) s(n-1) + b^n s(n) \equiv b^{n^2} A(b^n) \mod B(b^n).$$
As both terms in the left-hand side are positive, we do not need any correction term. For sufficiently large $b$, 
$$0 < (-\alpha_2) s(n-1) + b^n s(n) < B(b^n),$$
and also
$$b(n) s(n) > 0, $$
$$0 < (-\alpha_2) s(n-1) < b^n.$$
In conclusion if $\alpha_2 < 0$ the following representation works:
$$\forall \, n\geq n_0\,\,\,\, s(n) = \left \lfloor 
    \frac{ b^{n^2} A(b^n)  \bmod B(b^n)}{b^{n}} \right \rfloor.$$ 

 For the quotient -- quotient representation, consider the identity
 $$ B(X) S_{n}(X) = C_{n}(X) + X^{n+1} A(X)$$ 
and observe that 
$$C_{n}(X) = \alpha_2 s(n) - X s(n + 1),$$
again without intermediary terms.  It follows that if $\alpha_2 < 0$ then 
$$(-\alpha_2) s(n) + X s(n + 1) \equiv X^{n+1} A(X) \mod B(X),$$
$$(-\alpha_2) s(n) + b^n s(n + 1) \equiv b^{n(n+1)} A(b^n) \mod B(b^n).$$
As both terms in the left-hand side are positive, we do not need any correction term. For sufficiently large $b$, 
$$0 < (-\alpha_2) s(n) + b^n s(n + 1) < B(b^n),$$
and also
$$b^n s(n+1) > 0,$$
$$0 < (-\alpha_2) s(n) < b^n.$$
In conclusion if $\alpha_2 < 0$ the following representation works:
$$\forall \, n\geq n_0\,\,\,\, s(n) = \frac{1}{|\alpha_{2} |} \,\,\,\, ( b^{n(n+1)}A(b^n) ) \bmod B(b^n)       \bmod b^n.$$ 

\subsection{Degree 3, positive terms} 

In this subsection we mainly deduce the general representation for degree $3$ sequences. We also look for some particular case representations, but not in such exhaustive way as for the sequences of degree $2$.

Let $s \in \mathbb N ^ {\mathbb N}$ be a sequence with initial terms $s(0)$,  $s(1)$ and $s(2)$ and with recurrence polynomial $B(X) = \alpha_3 + \alpha_2 X + \alpha_1 X^2 +  X^3 \in \mathbb Z[X]$. As the polynomial $S_n(X)$ ends in
$$ \dots + s(2) X^{n-2} + s(1) X^{n-1} + s(0) X^n,$$
we conclude that $$A(X) = s(0) X^2 + (\alpha_1 s(0) + s(1)) X + (\alpha_2 s(0) + \alpha_1 s(1) + s(2)).  $$ According to Theorems \ref{Theoremmoddivrep} and \ref{Theoremmodmodrep}, the sequence $s$ has the following representations: 

 $$\forall \, n\geq n_0\,\,\,\, s(n) = \left \lfloor 
    \frac{\left ( b^{n + \left \lceil n/2 \right \rceil} + b^{n^2}A(b^n ) \right ) \bmod B(b^n)}{b^{2n}} \right \rfloor.$$

$$\forall \, n\geq n_0\,\,\,\, s(n) = \frac{1}{|\alpha_{3} |}  \left (b^{2n + \left \lceil n/2 \right \rceil} -  b^{n(n+1)} \sgn(\alpha_3) A(b^n) \right )\bmod B(b^n)       \bmod b^n.$$

 For the remainder -- quotient representation, consider the identity
 $$ B(X) S_{n-1}(X) = C_{n-1}(X) + X^{n} A(X)$$ 
and observe that 
$$C_{n-1}(X) = \alpha_3 s(n-1) + (\alpha_2 s(n-1) + \alpha_3 s(n-2)) X - X^2 s(n),$$
so that
$$X^2 s(n) - (\alpha_2 s(n-1) + \alpha_3 s(n-2)) X - \alpha_3 s(n-1) \equiv X^n A(X) \mod B(X),  $$
$$b^{2n} s(n) - (\alpha_2 s(n-1) + \alpha_3 s(n-2)) b^n - \alpha_3 s(n-1) \equiv b^{n^2} A(b^n) \mod B(b^n),  $$
When $\alpha_2 \leq 0$ and $\alpha_3 < 0$, both terms in the left-hand side are positive and we do not need any correction term. So in this special case, the following representation works: 

$$\forall \, n\geq n_0\,\,\,\, s(n) = \left \lfloor 
    \frac{ b^{n^2}A(b^n ) \bmod B(b^n)}{b^{2n}} \right \rfloor$$
for sufficiently large $b$. 

 For the quotient -- quotient representation, consider the identity
 $$ B(X) S_{n}(X) = C_{n}(X) + X^{n+1} A(X)$$ 
and observe that 

$$C_{n}(X) = \alpha_3 s(n) + (\alpha_2 s(n) + \alpha_3 s(n-1)) X - X^2 s(n +1),$$
so that 
$$|\alpha_3| s(n) + \sgn(\alpha_3) (\alpha_2 s(n) + \alpha_3 s(n-1)) X -  \sgn(\alpha_3) s(n+1) X^2 \equiv $$ $$ \equiv - \sgn(\alpha_3) X^{n+1}A(X) \mod B(X),$$ 
$$|\alpha_3| s(n) + \sgn(\alpha_3) (\alpha_2 s(n) + \alpha_3 s(n-1)) b^n -  \sgn(\alpha_3) s(n+1) b^{2n} \equiv $$ $$ \equiv - \sgn(\alpha_3) b^{n(n+1)}A(b^n) \mod B(b^n).$$
Again, if $\alpha_3 < 0$ and $\alpha_2 \leq 0$, all the terms in the left-hand side are positive and we don't need any correction term. In this particular case, the following representation
$$\forall \, n\geq n_0\,\,\,\, s(n) = \frac{1}{|\alpha_{3} |}  \,\,\,\,\,\, b^{n(n+1)}  A(b^n) \bmod B(b^n)       \bmod b^n$$
works for sufficiently large $b$.

\subsection{Degree 2, general case} 

If $s \in \mathbb Z^{\mathbb N}$ is given by the recurrence polynomial $B(X) = \alpha_2 + \alpha_1 X + X^2 \in \mathbb Z[X]$ with initial terms $s(0), s(1) \in \mathbb Z$, one has to find a constant $c \in \mathbb N$ and a sequence $t \in \mathbb N^{\mathbb N}$ such that for all $n \in \mathbb N$, $t(n) = s(n) + c^{n+1} $. 
As we know that a constant $c$ such that $|s(n)| < c^{n+1}$ always exists, and there are sufficient conditions for $c$, which are always satisfied by sufficiently large $c \in \mathbb N$, suppose that such a $c$ has been found and fixed. The sequence $t(n)$ has a recurrence given by the polynomial:
$$B_t (X) = (X^2 + \alpha_1X + \alpha_2)(X-c) = X^3 + (\alpha_1 - c) X^2 + (\alpha_2 - c \alpha_1) X - c\alpha_2,$$
with initial terms $t(0) = s(0) + c$, $t(1) = s(1) + c^2$ and $t(2) = s(2) + c^3$. 
We write 
$$B_t(X) = X^3 + \beta_1 X^2 + \beta_2 X + \beta_3 X$$
where $\beta_1 = \alpha_1 - c$, $\beta_2 = \alpha_2 - c\alpha_1$ and $\beta_3 = -c \alpha_2$. 

To represent the sequence $t$ with any of the methods remainder - quotient or remainder - remainder, we need the initial terms $t(0) = s(0) + c$, $t(1) = s(1) + c^2$ and $t(2) = s(2) + c^3 = - \alpha_1 s(1) - \alpha_2 s(2) + c^3$. We also need the degree $2$ polynomial $A(X)$ where: 
 $$A(X) = t(0) X^2 + (\beta_1 t(0) + t(1)) X + (\beta_2 t(0) + \beta_1 t(1) + t(2)).  $$
 We recall and compute:
 $$t(0) = s(0) + c,$$
 $$\beta_1 t(0) + t(1) = \alpha_1 (s(0) + c) + (s(1) - s(0) c ),$$
$$\beta_2 t(0) + \beta_1 t(1) + t(2) = \alpha_2(s(0) + c) + \alpha_1(s(1) - s(0) c) + (s(2) - s(1)c).  $$
Now we have all elements to represent the sequence $t$ and, consequently, to represent the sequence $s$. 

\section{Pell's equation}

Consider a non-square positive integer $ k $.

The Diophantine equation $$ X^2 - k Y^2 = 1 $$ is known as \textbf{Pell's equation} (see Barbeau \cite[Preface]{Barbeau}).
It is known that the pair $(x(0), y(0)) = (0,1)$ is always a solution. Also, there is a smallest non-trivial solution $(x(1), y(1))$ which is fundamental in the following sense: 
 If $ n $ is a positive integer, then the numbers $x(n)$ and $y(n)$ 
 defined by the relation: $$ x(n) \pm y(n) \sqrt{k} = {\left ( x(1) \pm y(1) \sqrt{k} \right ) }^n  $$ 
 build a solution of the Pell equation. Moreover, these are all solutions in the set of natural numbers.
In  Prunescu and Sauras-Atltuzarra \cite{PrunescuSaurasCRec1} it is proved that both sequences $x(n)$ and $y(n)$ satisfy the relation of recurrence: $$ s ( n + 2 ) = 2 x ( 1 ) s ( n + 1 ) - s ( n ) . $$ 
So the sequence $x(n)$ is determined by the relation above and the conditions $x(0) = 1$ and the given value of  $x(1)$, while similarly $y(n)$ is determined by the same relation and the information that $y(0) = 0$ together with the value of $y(1)$. 

In conclusion, to represent the sequences $x(n)$ and $y(n)$ we need the following polynomials. The recurrence polynomial for both sequences is $B(X) = X^2 - 2x(1) X + 1$. For the sequence $x(n)$, 
$A(X)  = X - x(1) $. For the sequence $y(n)$, $A(x) = y(1)$.

\section{Examples} 

\subsection{Degree 2, positive elements, negative free coefficient} 

The parameters of these sequences are $\alpha_1, \alpha_2 \in \mathbb Z$ and $s(0), s(1) \in \mathbb N$. Recall that:
$$B(X) = X^2 + \alpha_1 X + \alpha_2,$$
$$A(X) = s(0) X + (\alpha_1 s(0) + s(1)).  $$ 

We write down the special representations because they prove to work for all $n \geq 1$ if the chosen $b$ is sufficiently large. Not all representations will be given explicitly. 

We start with the sequence of \textbf{Fibonacci numbers} \href{https://oeis.org/A000045}{\texttt{OEIS A000045}}. The sequence is given by $\alpha_1 = \alpha_2 = -1$ and $s(0) = 0$, $s(1) = 1$.

\begin{proposition}\label{PropFormulaFibonacci} If $ n \geq 1 $, then $$ s(n) = \left \lfloor \frac{3^{n^2} \bmod (3^{2n} - 3^n - 1)}{3^n} \right \rfloor, $$ 
$$s(n) = 3^{n^2 + n} \bmod (3^{2n} - 3^n - 1) \bmod 3^n.$$
\end{proposition} 

We observe that also the following representation holds for the Fibonacci sequence for $n \geq 1$, and it follows also from the  presented theory:
$$s(n-1) = 2^{n^2} \bmod (2^{2n} - 2^n -1) \bmod 2^n. $$

The sequence of \textbf{Lucas numbers}  \href{https://oeis.org/A000032}{\texttt{OEIS A000032}} is given by $\alpha_1 = \alpha_2 = -1$ and $s(0) = 2$, $s(1) = 1$. 

\begin{proposition}\label{PropFormulaLucas} If $ n \geq 1 $, then  $$ s(n) = \left \lfloor \frac{(2 \cdot 4^{n^2 + n} - 4^{n^2} ) \bmod (4^{2n} - 4^n - 1)}{4^n} \right \rfloor, $$ 
$$s(n) = ( 2 \cdot 5^{n^2 + 2n} - 5^{n^2 + n}) \bmod (5^{2n} - 5^n - 1) \bmod 5^n.$$ \end{proposition}

The sequence of \textbf{Pell numbers} \href{https://oeis.org/A000129}{\texttt{OEIS A000129}} is given by 
$\alpha_1 = -2$, $\alpha_2 = -1$ and $s(0) = 0$, $s(1) = 1$.

\begin{proposition}\label{PropFormulaPellNumbers} If $ n \geq 1 $, then $$ s(n) = \left \lfloor \frac{ 4^{n^2} \bmod (4^{2n} - 2 \cdot 4^n - 1)}{4^n} \right \rfloor, $$ 
$$s(n) =  3^{n^2 + n} \bmod (3^{2n} - 2 \cdot 3^n - 1) \bmod 3^n.$$ \end{proposition}

The sequence of \textbf{Pell-Lucas numbers}  \href{https://oeis.org/A002203}{\texttt{OEIS A002203}} is given by $\alpha_1 = -2$, $\alpha_2 = -1$ and $s(0) = s(1) = 2$.

\begin{proposition} If $ n \geq 1$, then $$ s(n) = \left \lfloor \frac{(2 \cdot 5^{n^2 + n} - 2 \cdot 5^{n^2} ) \bmod (5^{2n} - 2 \cdot 5^n - 1)}{5^n} \right \rfloor, $$ 
$$s(n) = ( 2 \cdot 9^{n^2 + 2n} - 2 \cdot 9^{n^2 + n}) \bmod (9^{2n} - 2 \cdot 9^n - 1) \bmod 9^n.$$
\end{proposition}

\subsection{Degree 2, positive elements, positive free coefficient} 

We will see that in the case of the remainder - remainder representations, the quantities inside the first pairs of parentheses are always negative. The remainder of any integer division is always a natural number, and this is here of crucial importance.

Another remark: the representation Theorems work with both exponents $\lfloor n/2 \rfloor$ and $\lceil n/2 \rceil$. We will creatively use this freedom to build formulas which are true starting with a value $n_0$ which should be as low as possible. 

The sequence of all the non-negative integers  \href{https://oeis.org/A001477}{\texttt{OEIS A001477}} is given by $\alpha_1 = -2$, $\alpha_2 = 1$, $s(0) = 0$ and $s(1) = 1$. 

\begin{proposition}\label{PropFormulaNaturals} If $ n \geq 1$, then $$ \left \lfloor \frac{(4^{\lfloor n/2 \rfloor} + 4^{n^2}) \bmod (4^{2n} - 2 \cdot 4^n + 1)}{4^n} \right \rfloor = n . $$ 
If $n \geq 2$ then
$$(2^{n + \lceil n/2 \rceil} - 2^{n^2 + n} ) \bmod (2^{2n} - 2 \cdot 2^n + 1) \bmod 2^n = n.$$
\end{proposition}

The all-twos sequence  \href{https://oeis.org/A007395}{\texttt{OEIS A007395}} is given by $\alpha_1 = -2$, $\alpha_2 = 1$, $s(0) = 2$, $s(1) = 2$. 

\begin{proposition} If $ n \geq 1$, then $$  \left \lfloor \frac{(5^{\lceil n/2 \rceil} + 2 \cdot 5^{n^2 + n} - 2 \cdot 5^{n^2}) \bmod (5^{2n} - 2 \cdot 5^n + 1) }{5^n} \right \rfloor  = 2 .  $$ 
If $n \geq 2$, then
$$( 2^{n + \lfloor n/2 \rfloor} - 2 \cdot 2^{n^2 + 2n} + 2 \cdot 2^{n^2 + n} ) \bmod ( 2^{2n} - 2  \cdot 2^n + 1  ) \bmod 2^n = 2. $$
\end{proposition}

The sequence of \textbf{Mersenne numbers}  \href{https://oeis.org/A000225}{\texttt{OEIS A000225}} is given by $\alpha_1 = -3$, $\alpha_2 = 2$, $s(0) = 0$, $s(1) = 1$. 

\begin{proposition}\label{PropFormulaMersenne} If $ n \geq 1 $, then $$\left \lfloor \frac{(6^{\lfloor n/2 \rfloor} + 6^{n^2}) \bmod (6^{2n} - 3 \cdot 6^n + 2)}{6^n} \right \rfloor = 2^n - 1 . $$ 
If $n \geq 2$, then
$$ \frac{1}{2} \cdot \left ( (4 ^ {n +  \lceil n / 2 \rceil } - 4 ^ {n ^ 2 + n} ) \bmod (4 ^ {2n} - 3 \cdot 4 ^ n + 2 ) \bmod 4 ^ n \right ) = 2^n - 1.$$
\end{proposition}

The sequence that maps each non-negative integer $ n $ into $ 2^n + 1 $ \href{https://oeis.org/A000051}{\texttt{OEIS A000051}} is given by $\alpha_1 = -3$, $\alpha_2 = 2$, $s(0) = 2$, $s(1) = 3$. 

\begin{proposition} If $ n \geq 1 $, then $$\left \lfloor \frac{( 7^{\lceil n/2 \rceil} +  2 \cdot 7^{n^2 + n} - 3 \cdot 7^{n^2}) \bmod (7^{2n} - 3 \cdot 7^n + 2)}{7^n} \right \rfloor  = 2^n + 1 , $$ 
$$ \frac{1}{2} \cdot  \left ( ( 7^{n + \lceil n/2 \rceil } - 2 \cdot  7^{n^2 + 2n} + 3 \cdot  7^{n^2 + n} ) \bmod ( 7^{2n}  -3 \cdot  7^n + 2 ) \bmod 7^n \right ) = 2^n + 1.$$
\end{proposition}

Consider Pell's equation:
$$  X^2 - 7 Y^2 = 1. $$

It is easy to see that the smallest non-trivial solution is the pair $ ( x(1) , y(1) ) = ( 8 , 3 ) $.  The sequence $x(n)$ is given by $\alpha_1 = -16$, $\alpha_2 = 1$, $x(0) = 1$ and $x(1) = 8$, and is indexed as \href{https://oeis.org/A001081}{\texttt{OEIS A001081}}. The sequence $y(n)$ is given by $\alpha_1 = -16$, $\alpha_2 = 1$, $y(0) = 0$ and $y(1) = 3$, and is indexed as \href{https://oeis.org/A001080}{\texttt{OEIS A001080}}. For the following representations I did not look for the smallest value of $b$. I just observed that $b = 256$ does it.

\begin{proposition} If $n \geq 0$, then
$$ x(n) = \left \lfloor  \frac{( 256^{\lfloor n / 2 \rfloor} +  256^{n^2 + n}  - 8 \cdot 256^{n^2})  \bmod (256^{2n} -16 \cdot 256^n + 1 )} {256^n} \right \rfloor, $$

$$ x(n) = (256 ^ {n + \lceil n / 2 \rceil} -  256 ^ {n ^ 2 + 2n} +8 \cdot 256 ^ {n ^ 2 + n} ) \bmod(256 ^ {2n} -16 \cdot 256 ^ n + 1 ) \bmod 256 ^ n, $$ 

$$ y(n) = \left \lfloor  \frac{( 256^{\lfloor n / 2 \rfloor} + 3 \cdot  256^{n^2} ) \bmod (256^{2n} -16 \cdot 256^n + 1 )} {256^n} \right \rfloor, $$

$$ y(n) = (256 ^ {n + \lceil n / 2 \rceil} - 3 \cdot 256 ^ {n ^ 2 + n}  ) \bmod(256 ^ {2n} -16 \cdot 256 ^ n + 1 ) \bmod 256 ^ n, $$ 
    
\end{proposition}

\subsection{Degree 2, integers}

In this subsection, we obtain formulas for two Lucas sequences that take both positive and negative values.

The first sequence is the sequence of \textbf{generalized Gaussian Fibonacci integers}  \href{https://oeis.org/A088137}{\texttt{OEIS A088137}}. It is given by $\alpha_1 = -2$, $\alpha_2 = 3$, $s(0) = 0$, $s(1) = 1$. For $c = 3$, we define the sequence $t(n) = s(n) + c^{n+1}$ and we have $t(0) = 3$, $t(1) = 10 $ and $t(2) = 29$. The recurrence of the sequence $t(n)$ is given by the polynomial:
$$B(X) = (X^2 - 2X + 3)(X-3) = X^3 - 5X^2 + 9X - 9.$$
The corresponding denominator of $t(n)$ is:
$$A(X) = 3X^2 - 5X + 6.$$
As $\beta_2 = 9 > 0$, we need a correction term. 

\begin{proposition} If $ n \geq 1 $, then $t(n) =$
$$ = \left \lfloor \frac{(32^{n + \lceil n/2 \rceil } + 3 \cdot 32^{n^2 + 2n} - 5 \cdot 32^{n^2 + n} + 6 \cdot 32^{n^2}) \bmod (32^{3n} - 5\cdot 32^{2n} + 9 \cdot 32^n - 9)}{32^{2n}} \right \rfloor, $$
$$s(n) = t(n) - 3^{n+1},$$
$$s(n) = \frac{1}{9} \cdot \big ( (128^{2n + \lfloor n/2} + 3 \cdot 128^{n^2 + 3n} - 5 \cdot 128^{n^2 + 2n} + 6 \cdot 128^{n^2 + n} ) \bmod $$ $$  \bmod (128^{3n} - 5 \cdot 128^{2n} + 9 128^n - 9) \bmod 128^n \big ) - 3^{n+1}.$$
\end{proposition}

The second sequence is \href{https://oeis.org/A002249}{\texttt{OEIS A002249}}. It is defined by $\alpha_1 = -1$ and $\alpha_2 = 2$, with initial terms $s(0) = 2$ and $s(1) = 1$. For $c = 2$, we define the sequence $t(n) = s(n) + c^{n+1}$ and we have $t(0) = 4$, $t(1) = 5$, $t(2) = 5$.  The recurrence of the sequence $t(n)$ is given by the polynomial:
$$(X^2 - X + 2)(X-2) = X^3 -3X^2 +4X -4.$$
The corresponding denominator of $t(n)$ is:
$$A(X) = 4X^2 - 7X + 6.$$
As $\beta_2 = 4 > 0$, we need a correction term. 

\begin{proposition} If $ n \geq 1$, then $s(n) = $
$$ = \left \lfloor \frac{((8^{n + \lceil n/2 \rceil } + 4 \cdot 8^{n^2 + 2n} - 7 \cdot 8^{n^2 + n} + 6 \cdot 8^{n^2}) \bmod (8^{3n} - 3\cdot 8^{2n} + 4 \cdot 8^n - 4)}{8^{2n}} \right \rfloor  - 2^{n+1} , $$
$$s(n) = \frac{1}{4} \cdot \big ( (32^{2n + \lfloor n/2 \rfloor} 
+ 4\cdot 32^{n^2 + 3n} - 7 \cdot 32^{n^2 + 2n} + 6 \cdot 32^{n^2 + n} ) \bmod $$ $$\bmod (32^{3n} - 3 \cdot 32^{2n} + 4 \cdot 32^{n} - 4) \bmod 32^n \big ) - 2^{n+1}. $$
\end{proposition} 

\subsection{Degree 3, positive elements, negative coefficients}

We finally apply the theory to some C-recursive natural sequences of degree three, whose recursions do not contain positive coefficients. Consequently these representations do not need correction terms. 
 
The sequence of \textbf{Tribonacci numbers} \href{https://oeis.org/A000073}{\texttt{OEIS A000073}} is the sequence $ s $ such that $s(0) = s(1) = 0$, $s(2) = 1$ and $$ s(n) = s(n-1) + s(n-2) + s(n-3) $$ for every integer $ n \geq 3 $.

\begin{proposition} If  $ n \geq 1 $, then 
$$ s(n) =  \left \lfloor \frac{2^{n^2} \bmod (2^{3n} - 2^{2n} - 2^n - 1)}{2^{2n}} \right \rfloor, $$ 
$$s(n) = 2^{n^2 + n} \bmod (2^{3n} - 2^{2n} - 2^n - 1) \bmod 2^n. $$
\end{proposition}

The sequence of \textbf{Padovan numbers} \href{https://oeis.org/A000931}{\texttt{OEIS A000931}}  is  the sequence $ s $ such that $s(0) = 1$, $s(1) = s(2) = 0$ and $$ s(n) = s(n-2) + s(n-3) $$ for every integer $ n \geq 3 $.

\begin{proposition} If $ n \geq 1 $, then
$$ s(n) = \left \lfloor \frac{(2^{n^2 +2 n} - 2^{n^2}) \bmod (2^{3n} - 2^n - 1)}{2^{2n}} \right \rfloor , $$ 
$$s(n) = (2^{n^2 + 3n} - 2^{n^2 + n}) \bmod (2^{3n} - 2^n - 1) \bmod 2^n. $$
\end{proposition}

The \textbf{Narayana's cows sequence} \href{https://oeis.org/A000930}{\texttt{OEIS A000930}}  is  the sequence $ s $ such that $s(0) = s(1) = s(2) = 1$ and $$ s(n) = s(n-1) + s(n-3) $$ for every integer $ n \geq 3 $.

\begin{proposition} If $ n \geq 1 $, then
$$ s(n) = \left \lfloor \frac{3^{n^2 +2 n} \bmod (3^{3n} - 3^{2n} - 1)}{3^{2n}} \right \rfloor , $$ 
$$s(n) = 2^{n^2 + 3n}  \bmod (2^{3n} - 2^{2n} - 1) \bmod 2^n. $$
\end{proposition}

\section{Conclusions}

\begin{enumerate}
    \item After the quotient - remainder representation of the C-recursive integer sequences done in \cite{PrunescuSaurasCRec1}, here is shown that other two related representations exist: a remainder - quotient representation and a remainder - remainder representation. 
    \item Instead of computing the quotient of the polynomials $A(X)$ and $B(X)$, for the two new representations one considers the remainder of this division. 
    \item Differently from the quotient - remainder representation, in most cases one has to add a correction term to the polynomial related to the numerator of the generating function. 
    \item Every C-recursive sequence of integers, and specially the oscillating sequences, can be written down as the difference between such a modular term and a geometric progression.
    \item Among the applications we count the sequences of solutions $(x(n), (y(n))$ of a Pell equation, the Lucas sequences and various classical sequences of order $3$. 
\end{enumerate}


\begin{thebibliography}{99}

\bibitem{Barbeau}
\newblock E.\ J.\ Barbeau,
\newblock \textit{Pell's Equation},
\newblock Problem Books in Mathematics,
\newblock Springer, New York,
\newblock 2002.

\bibitem{BicknellJohnsonHoggatt}
\newblock M.\ Bicknell-Johnson,  V.\ E.\ Hoggatt,
\newblock Fibonacci convolution sequences,
\newblock \textit{Fibonacci Quarterly} \textbf{15} (1977).

\bibitem{BorweinCrandall}
\newblock J.\ M.\ Borwein,  R.\ E.\ Crandall,
\newblock Closed Forms: What They Are and Why We Care,
\newblock \textit{Notices of the American Mathematical Society} \textbf{60} (2013).
\newblock \url{https://doi.org/10.1090/NOTI936}.

\bibitem{Enderton}
\newblock H.\ B.\ Enderton,
\newblock \textit{A Mathematical Introduction to Logic} (2nd ed.),
\newblock Academic Press,
\newblock 2001.

\bibitem{Grigorieva}
\newblock E.\ Grigorieva,
\newblock \textit{Methods of Solving Number Theory Problems},
\newblock Birkhäuser,
\newblock 2018.

\bibitem{Guy}
\newblock R.\ K.\ Guy,
\newblock \textit{Unsolved Problems in Number Theory} (3rd ed.),
\newblock Springer, New York,
\newblock 2004.
\newblock \url{https://doi.org/10.1007/978-0-387-26677-0} 

\bibitem{Hankin}
\newblock Paul Hankin, Far\'e Rideau,
\newblock \textit{A Novel and Efficient Way to Compute Fibonacci Numbers,}
\newblock 2018.
\newblock Blog post.
\newblock \url{https://blog.paulhankin.net/fibonacci2}

\bibitem{KrizekEtAl}
\newblock M.\ K\v{r}\'{i}\v{z}ek, F.\ Luca,  L.\ Somer,
\newblock \textit{17 Lectures on Fermat Numbers: from Number Theory to Geometry},
\newblock CMS Books in Mathematics,
\newblock Springer, New York,
\newblock 2001.

\bibitem{Marchenkov}
\newblock S.\ S.\ Marchenkov,
\newblock Superpositions of Elementary Arithmetic Functions,
\newblock \textit{Journal of Applied and Industrial Mathematics} \textbf{3} (2007).

\bibitem{Mazzanti}
\newblock S.\ Mazzanti,
\newblock Plain Bases for Classes of Primitive Recursive Functions,
\newblock \textit{Mathematical Logic Quarterly} \textbf{48} (2002).

\bibitem{Mendelson}
\newblock E.\ Mendelson,
\newblock \textit{Introduction to Mathematical Logic} (6th ed.),
\newblock Taylor \& Francis,
\newblock 2015.

\bibitem{Oitavem}
\newblock I.\ Oitavem,
\newblock New recursive characterizations of the elementary functions and the functions computable in polynomial space,
\newblock \textit{Revista Matemática de la Universidad Complutense de Madrid} \textbf{10} (1997).

\bibitem{PetkovsekEtAl}
\newblock M.\ Petkovšek, H.\ S.\ Wilf and D.\ Zeilberger,
\newblock \textit{A = B},
\newblock A.\ K.\ Peters / CRC Press,
\newblock 1996.
\newblock \url{http://www.math.upenn.edu/~wilf/AeqB.html}

\bibitem{PetkovsekZakrajsek}
\newblock M.\ Petkovšek, H.\ Zakrajšek,
\newblock Solving linear recurrence equations with polynomial coefficients,
\newblock in J.\ Blümlein and C.\ Schneider, eds.,
\newblock \textit{Computer algebra in quantum field theory. Integration, summation and special functions.},
\newblock Springer,
\newblock 2013,
\newblock pp.\ 259--284.

\bibitem{PrunescuSaurasCRec1} 
\newblock M.\ Prunescu, L.\ Sauras-Altuzarra,
\newblock \textit{On the representation of C-recursive integer sequences by arithmetic terms}
\newblock 	arXiv:2405.04083 
\newblock \url{https://arxiv.org/abs/2405.04083}.

\bibitem{PrunescuSaurasAltuzarra}
\newblock M.\ Prunescu, L.\ Sauras-Altuzarra,
\newblock \textit{An arithmetic term for the factorial function},
\newblock Examples \& Counterexamples, \textbf{5} (2024).
\newblock \url{https://doi.org/10.1016/j.exco.2024.100136}.

\bibitem{Rosen}
\newblock K.\ H.\ Rosen,
\newblock \textit{Elementary Number Theory and Its Applications} (6th ed.),
\newblock Addison-Wesley,
\newblock 2011.

\bibitem{SaurasAltuzarra}
\newblock L.\ Sauras-Altuzarra,
\newblock \textit{Hypergeometric closed forms},
\newblock master thesis,
\newblock Vienna University of Technology,
\newblock 2018.
\newblock \url{https://doi.org/10.25365/thesis.57260} 

\bibitem{Shunia}
\newblock Joseph M. Shunia,
\newblock \textit{Polynomial Quotient Rings and Kronecker Substitution for Deriving Combinatorial Identities.}
\newblock  arXiv:2404.00332
\newblock \url{https://arxiv.org/abs/2404.00332}

\bibitem{Stanley}
\newblock R.\ P.\ Stanley,
\newblock \textit{Enumerative Combinatorics} (vol.\ 1, 2nd ed.),
\newblock Cambridge University Press,
\newblock 2011.

\bibitem{Encyclopedia}
\newblock Various,
\newblock ``Lucas sequence'',
\newblock Encyclopedia of Mathematics.
\newblock \url{https://encyclopediaofmath.org/wiki/Lucas_sequence}

\bibitem{VereschchaginShen}
\newblock N.\ K.\ Vereschchagin,  A.\ Shen,
\newblock \textit{Computable Functions} (translated by V.\ N.\ Dubrovskii),
\newblock American Mathematical Society,
\newblock 2002.

\bibitem{Weisstein}
\newblock E.\ W.\ Weisstein,
\newblock ``Eventually'',
\newblock from MathWorld -- A Wolfram Web Resource.
\newblock \url{https://mathworld.wolfram.com/Eventually.html}

\bibitem{Weisstein2}
\newblock E.\ W.\ Weisstein,
\newblock ``Generating Function'',
\newblock from MathWorld -- A Wolfram Web Resource.
\newblock \url{https://mathworld.wolfram.com/GeneratingFunction.html}

\bibitem{Wilf}
\newblock H.\ S.\ Wilf,
\newblock \textit{Generatingfunctionology},
\newblock A.\ K.\ Peters,
\newblock 2006.

\end{thebibliography}
\end{document}